\documentclass[pdflatex,sn-mathphys-num]{sn-jnl}

\usepackage{graphicx}%
\usepackage{multirow}%
\usepackage[all]{xy}
\usepackage{amsmath,amssymb,amsfonts}%
\usepackage{amsthm}%
\usepackage{mathrsfs}%
\usepackage[title]{appendix}%
\usepackage{xcolor}%
\usepackage{textcomp}%
\usepackage{manyfoot}%
\usepackage{booktabs}%
\usepackage{algorithm}%
\usepackage{algorithmicx}%
\usepackage{algpseudocode}%
\usepackage{listings}%
\usepackage{soul}
\usepackage{enumitem}

\theoremstyle{thmstyleone}%
\newtheorem{theorem}{Theorem}
\theoremstyle{thmstyletwo}%
\newtheorem{remark}{Remark}%
\newtheorem{lemma}[theorem]{Lemma}
\theoremstyle{thmstylethree}%
\newtheorem{definition}{Definition}%
\newtheorem{corollary}[theorem]{Corollary}
\begin{document}

\title[Topological pressures of a factor map for nonautonomous iterated function systems]{Topological pressures of a factor map for nonautonomous iterated function systems}


\author*[1]{\fnm{Yujun} \sur{Ju}}\email{yjju@ctbu.edu.cn}
\equalcont{These authors contributed equally to this work.}

\author[1]{\fnm{Lingbing} \sur{Yang}}\email{yanglingbing@ctbu.edu.cn}
\equalcont{These authors contributed equally to this work.}

\affil*[1]{\orgdiv{School of Mathematics and Statistics}, \orgname{Chongqing Technology and Business University}, \orgaddress{\city{Chongqing}, \postcode{400067}, \country{People's Republic of China}}}


\abstract{This study focuses on the topological pressure of nonautonomous iterated function systems defined on a compact metric space. We establish an inequality relating two topological pressures associated with a factor map of nonautonomous iterated function systems.}




\maketitle

\section{Introduction}\label{intro}

The concept of topological entropy, first introduced by Adler et al. \cite{AKM}, serves as a fundamental invariant in dynamical systems, quantifying the exponential growth rate of distinguishable orbits over time to reveal the system's underlying complexity. Bowen \cite{Bowen1} considered a surjective continuous map $\pi: (X, T) \rightarrow (Y, S)$, establishing the inequality
\[
h(T) \leq h(S)+\sup _{y \in Y} h\left(\pi^{-1}(y),T\right),
\]
where $h(K, T)$ denotes the topological entropy of a compact subset $K$ of $X$. Kolyada and Snoha \cite{KS} extended this result in 1996 to non-autonomous discrete dynamical systems. For equicontinuous non-autonomous systems $(X,f_{1, \infty})$ and $(Y,g_{1, \infty})$, they proved that if $\{\pi_i\}_{i=1}^{\infty}$ is a semiconjugacy and $\pi_i \in \{\phi_1, \phi_2, \ldots, \phi_k\}, \, k < \infty$, for every $i \geq 1$, then
\[
h(f_{1,\infty})\leq h(g_{1,\infty})+\max_{i}\sup_{y\in Y}H(f_{1,\infty};\pi_{i}^{-1}(y)),
\]
where $H(f_{1, \infty}; \pi_{i}^{-1}(y))$ denotes the topological sup-entropy of the sequence $f_{1, \infty}$ on the set $\pi_{i}^{-1}(y)$.

Topological pressure, a generalization of topological entropy introduced by Ruelle \cite{Ruelle1973}, measures orbit complexity of an iterated map on a given potential function and has since become a central concept in thermodynamic formalism and multifractal analysis. Later, Walters \cite{Walters1975} extended this concept to compact spaces with continuous transformations. More recently, Lin et al. \cite{Lin2018} formulated definitions of topological pressures for free semigroup actions, drawing on by Bufetov's \cite{Bufe} definition of topological entropy for such actions, and established a partial variational principle. Building on Bi{\'s}'s entropy for free semigroup actions \cite{Bis}, Zhao and Yang \cite{Zhao} established an analogue of Bowen's inequality for topological pressures with factor maps.

In 2019, Ghane and Nazarian \cite{ghane2019} introduced topological entropy and topological pressure for non-autonomous iterated function systems (or NAIFSs for short), a framework that generalizes both the concepts of free semigroup actions and non-autonomous discrete dynamical systems. Unlike usual (autonomous) iterated function systems, NAIFSs evolve over time, which is essential in various applications. Generalized Cantor sets that studied by Robinson and Sharples \cite{Robinson2012} serve as examples of attractors of NAIFSs. Rempe-Gillen and Urba\'nski \cite{Rempe2016} investigated the Hausdorff dimension of limit sets of NAIFSs, demonstrating that under a suitable restriction on the growth of the number of contractions used at each step, it is  determined by an equation known as Bowen's formula. They also extended this formula to a class of infinite alphabet systems, addressed Hausdorff measures for finite systems, and examined the continuity of topological pressure and Hausdorff dimension for both finite and infinite systems.

Despite these developments, the behavior of topological pressure with factor maps for NAIFSs remains insufficiently explored. In particular, it is not yet well understood how classical results--such as Bowen's inequality and its extension by Zhao and Yang--can be formulated or generalized within this broader, time-dependent framework. Bridging this gap is essential for advancing the thermodynamic formalism of NAIFSs and deepening our understanding of the relationship between system structure and dynamical complexity in non-autonomous settings.

In this paper, we study foundational properties of topological pressure for NAIFSs, including the power rule and monotonicity. The principal contribution of this work is a novel inequality between two topological pressures with a factor map in the NAIFS setting. This result is derived by introducing the notion of topological sup-entropy for NAIFSs, which plays a central role in the derivation. Rather than being a mere generalization of non-autonomous discrete dynamical systems and free semigroup actions, NAIFSs serve as a natural bridge between them. By adopting and extending techniques originally developed in the context of non-autonomous discrete dynamical systems, we derive an inequality for the topological pressure with factor maps of NAIFSs, which in turn yields a more refined approximation of the corresponding inequality for free semigroup actions, thereby improving the result of Zhao and Yang \cite{Zhao}.

The structure of this paper is as follows.  In Sect. 2, we review relevant preliminaries. In Sect. 3, we discuss fundamental properties of topological pressure for NAIFSs. In Sect. 4, we provide an inequality relating two topological pressures and derive a formula involving factor maps of NAIFSs.

\section{Preliminaries}

Let $\mathbb{Z}$ denote the set of integers, $\mathbb{N}$ the set of positive integers, and $\mathbb{R}^{+}$ the set of non-negative real numbers. Let $X$ be a set and let $\Phi$ consist of a sequence $\left\{\Phi^{(j)}\right\}_{j \geq 1}$ of collections of maps, where each $\Phi^{(j)}=\left\{f_{i}^{(j)}: X \rightarrow X\right\}_{i \in I^{(j)}}$, with $I^{(j)}$ being a nonempty finite index set. The pair $(X, \Phi)$ is a non-autonomous iterated function system (or NAIFS for short). To simplify notation, define the following symbolic spaces for $m \in \mathbb{Z}$ and $n \in \mathbb{N}$:
\[
I^{m, n}:=\prod_{j=0}^{n-1} I^{(m+j)},
\quad I^{m, \infty}:=\prod_{j=m}^{\infty} I^{(j)},
\quad I^{-\infty, \infty}:=\prod_{j=-\infty}^{\infty} I^{(j)}.
\]
Denote by $\Sigma_m$ the set of all two-sided infinite sequences of symbols $ 0, \ldots, m-1$.
It is clear that $ I^{-\infty, \infty}=\Sigma_m $ if $I^{(j)}=\{0, \ldots, m-1\}$ for all $j \in \mathbb{Z}.$
If $w \in I^{m, n}$, then $w$ is called a finite word and its length is $n$ and denoted by $|w|:=n .$ While, each word $\omega \in I^{m, \infty}$ or $I^{-\infty, \infty}$ is called an infinite word and its length is infinity and denoted by $|\omega|:=\infty.$ 
Let $\omega=(\ldots,\omega_{-1},\omega_{0},\omega_{1},\ldots) \in I^{-\infty,\infty}, w \in I^{m,n},$  write $\left.\omega\right|_{[m, m+n-1]}=w$ if $w=\omega_{m} \omega_{m+1} \cdots \omega_{m+n-1}$.

Let $ i \geq m$ and $k \in \mathbb{N} $. For a finite word $ w=w_{m} w_{m+1} \cdots w_{m+n-1}\in I^{m, n} $
with $i+k \leq m+n$, or for an infinite word $w \in I^{m,\infty}(I^{-\infty,\infty})$, we define $f_{w}^{m, 0}:=id_{X}$, where $id_{X}$ is the identity map on $X$. Set
\[
f_{w}^{i, k}:=f_{w_{i+k-1}}^{(i+k-1)} \circ \cdots \circ f_{w_{i+1}}^{(i+1)} \circ f_{w_{i}}^{(i)}\quad \text {and} \quad f_{w}^{i, -k}:=(f_{w}^{i, k})^{-1}.
\]

Let an $\mathrm{NAIFS}$ $(X, \Phi)$ and $n \geq 1$ be given. 
Denote by $\left(X, \Phi^{n}\right)$ the NAIFS defined by the sequence $\left\{\Phi^{(j, n)}\right\}_{j \geq 1}$, where $\Phi^{(j, n)}$ is the collection $\left\{f_{w_{j}^{*}}^{(j, n)}\right\}_{w_{j}^{*} \in I^{(j, n)}}$, $I^{(j, n)}:=\left\{w_{j}^{*} \in I^{(j-1) n+1, n}\right\}$. Note that $I^{(j, n)}=I^{(j-1) n+1, n}$ and 
\[
f_{w_{j}^{*}}^{(j, n)}:=f_{w_{j n}}^{(j n)} \circ \cdots \circ f_{w_{(j-1) n+2}}^{((j-1) n+2)} \circ f_{w_{(j-1) n+1}}^{((j-1) n+1)} \text { for } w_{j}^{*}=w_{(j-1) n+1} w_{(j-1) n+2} \cdots w_{j n}.
\] 
Take $I_{*}^{m, k}:=\prod_{j=0}^{k-1} I^{(m+j, n)},$ then $\#\left(I_{*}^{1, m}\right)=\#\left(I^{1, m n}\right),$ where $\#(A)$ is the cardinal number of the set $A .$ 
For any $w=w_{1} w_{2} \cdots w_{m n} \in I^{1, m n}$ and $1 \leq j \leq m,$ denote $w_{(j-1) n+1} w_{(j-1) n+2} \cdots w_{j n}$ by $w_{j}^{*} \in I^{(j, n)},$ then $w=w_{1}^{*} w_{2}^{*} \cdots w_{m}^{*} \in I_{*}^{1, m} .$ 
For simplicity, we denote elements in $I_{*}^{1, m}$ by $w^{*}$ and use analogous notation for other sequences of objects related to an NAIFS.

Let $(X, \Phi)$ be an NAIFS of continuous maps on a compact metric space $(X,d)$.
For finite (infinite) word $w=w_m w_{m+1} \cdots w_{m+n-1}\left(w=w_m w_{m+1} \cdots\right) \in$ $I^{m, n}\left(I^{m, \infty}\right)$ and $1 \leq k \leq|w|(1 \leq k<\infty)$ we introduce the Bowen metrics on $X$ 
\[
d_{w,k}(x, y):=\max _{0 \leq j \leq k} d\left(f_w^{1,j}(x), f_w^{1,j}(y)\right).
\]

Fix $w \in I^{1,n}$ for some $n \geq 1.$ A subset $E$ of the space $X$ is called $(n, w, \epsilon; \Phi)$-separated, if for any two distinct points $x,y \in E, d_{w,n}(x,y) > \epsilon$. Also, a subset $F$ of the space $X$, $(n, w, \epsilon; \Phi)$-spans another subset
$K \subseteq X$, if for each $x \in K$ there is a $y \in F$ such that $d_{w,n}(x,y) \leq \epsilon$. For a subset $Y $ of $X$ we define $s_n(Y; w, \epsilon, \Phi)$, as the maximal cardinality of an $(n, w, \epsilon; \Phi)$-separated set in $Y$ and $r_n(Y; w, \epsilon, \Phi) $ as the minimal cardinality of a set in $Y$ which $(n, w, \epsilon; \Phi)$-spans $Y$. If $Y = X$ we sometime suppress $Y$ and shortly write $s_n(w, \epsilon, \Phi)$ and $r_n(w, \epsilon, \Phi)$.
Let $C(X,\mathbb{R})$ be the space of real-valued continuous functions of $X.$ For $\psi\in C(X,\mathbb{R})$ and finite word $w\in I^{m,n}$ we denote $\Sigma_{j=0}^n\psi(f_w^{m,j})(x)$ by $S_{w,n}\psi(x).$ 
For $\epsilon>0, n \ge1, w\in I^{1,n}$ and $\psi\in C(X,\mathbb{R})$, put
\[
Q_n(\Phi;w,\psi,\epsilon):=\inf_F\left\{\sum_{x\in F}e^{S_{w,n}\psi(x)}:F\text{ is an }(n,w,\epsilon;\Phi)\text{-spanning set for } X\right\},
\]
\[
P_n(\Phi;w,\psi,\epsilon):=\sup_E\left\{\sum_{x\in E}e^{S_{w,n}\psi(x)}:E\text{ is an }(n,w,\epsilon;\Phi)\text{-separated set for }X\right\},
\]
and taking
\[
Q_n(\Phi;\psi,\epsilon):=\frac{1}{\#(I^{1,n})}\sum_{w\in I^{1,n}}Q_n(\Phi;w,\psi,\epsilon).
\]
\[
P_n(\Phi;\psi,\epsilon):=\frac{1}{\#(I^{1,n})}\sum_{w\in I^{1,n}}P_n(\Phi;w,\psi,\epsilon).
\]
\begin{definition}[\cite{ghane2019}]
For $\psi \in C(X,\mathbb{R})$, the topological pressure of an NAIFS $(X,\Phi)$ with respect to $\psi$ is defined as
\[
P(\Phi,\psi):=\lim_{\epsilon\to 0}Q(\Phi;\psi,\epsilon)=\lim_{\epsilon\to0}P(\Phi;\psi,\epsilon),
\]
where 
\[
Q(\Phi;\psi,\epsilon):=\lim\sup_{n\to\infty}\frac{1}{n}\log Q_n(\Phi;\psi,\epsilon),
\]
\[
P(\Phi;\psi,\epsilon):=\lim\sup_{n\to\infty}\frac1n\log P_n(\Phi;\psi,\epsilon).
\]
\end{definition}
\begin{remark}
(1) When $\#(I^{(j)}) = 1$ and $\Phi^{(j)}=\{f_1^{(j)}\}$ for every $j \in \mathbb{N}$, then we get the definition of topological pressure for non-autonomous discrete dynamical system $(X,f_{1,\infty})$ \cite{Huang2008,Kong2015}, where $f_{1,\infty}$ is the sequence $\{f_1^{(j)}\}_{j=1}^{\infty}$. Additionally, if $f_1^{(j)}=f$ for every $j\geq1$, this definition is just the topological pressure $P(f,\psi)$ for autonomous dynamical system $(X, f).$ 

(2) When $ \Phi^{(i)}=G:=\{f_{0},f_{1},\ldots,f_{m-1}\} $ for all $ i \in \mathbb{N}$, then we get the definition of topological pressure for semigroup action $ (X,G) $ with generator set $G$.

(3) When $\psi = 0$, $P( \Phi, \psi )$ reduces to the topological entropy of NAIFSs defined in  \cite{ghane2019} and is denoted by $h(\Phi).$ 
\end{remark}

\section{Dynamical properties}\label{section3}

\begin{theorem}\label{th1}
Let $(X,\Phi)$ be an NAIFS on a compact metric space $(X,d).$ If $\psi,\varphi\in C(X,\mathbb{R}),\epsilon>0$ and $c\in\mathbb{R}$, then the following are true.

(i) If $\varphi \leq \psi$, then $P(\Phi,\varphi)\leq P(\Phi,\psi)$ holds. In
particular,
\[
h(\Phi)+\inf \varphi \leq P(\Phi,\varphi)\leq h(\Phi)+ \sup \varphi.
\]

(ii) For any $\varphi\in C(X,\mathbb{R})$,
\[
P(\Phi,\varphi)
\begin{cases}=+\infty&\mathrm{if~}h(\Phi)=+\infty,\\<+\infty&\mathrm{if~}h(\Phi)<+\infty.\end{cases}
\]

(iii) 
\[
|P(\Phi;\varphi,\epsilon)-P(\Phi;\psi,\epsilon)|\leq\\\|\varphi-\psi\|,\]
and so if $P(\Phi,\cdot)$ is finite,
then 
\[
|P(\Phi,\varphi)-P(\Phi, \psi ) | \leq \| \varphi - \psi \|. \]
In other words,  $P(\Phi, \cdot)$ is a continuous function on $C(X,\mathbb{R}).$

(iv) $P(\Phi, \cdot , \epsilon )$ is convex, and so if
$P(\Phi,\cdot)<\infty$, then $P(\Phi,\cdot)$ is convex.

(v) $P(\Phi,\varphi+c)=P(\Phi,\varphi)+c.$

(vi)
\[
P(\Phi,\varphi+\psi)\leq P(\Phi,\varphi)+P(\Phi,\psi)+\limsup_{n \to \infty}\frac{\sum_{i=1}^{n}\log \#I^{(i)}}{n}.
\]

(vii)
\begin{align*}
&P(\Phi,c\varphi)\leq cP(\Phi,\varphi)+(c-1)\limsup_{n \to \infty}\frac{\sum_{i=1}^{n}\log \#I^{(i)}}{n}\quad ~if~c\geq 1\quad and~\\
&P(\Phi,c\varphi)\geq cP(\Phi,\varphi)+(c-1)\liminf_{n \to \infty}\frac{\sum_{i=1}^{n}\log \#I^{(i)}}{n}\quad ~if~ c\leq 1.
\end{align*}

(viii) 
\[
-2\liminf_{n \to \infty}\frac{\sum_{i=1}^{n}\log \#I^{(i)}}{n}-P\left(\Phi,\vert\varphi\vert\right)\leq P\left(\Phi,\varphi\right)\leq P\left(\Phi,\vert\varphi\vert\right).
\]
\begin{proof}
$(i)$ The result follows directly from the definition of topological pressure.

$(ii)$ According to $(i)$, we know $\inf\varphi+h(\Phi)\leq P(\Phi,\varphi) \leq h(\Phi)+\sup\varphi.$ Note that $\inf\varphi,\sup\varphi$ are finite since $X$ is compact and $\varphi\in C(X,\mathbb{R}).$ Thus, $P(\Phi,\varphi)=\infty$ if and only if $h(\Phi)=\infty.$

$(iii)$ By the inequality $\sup_{E}ab\leq\sup_{E}a\sup_{E}b,a,b\geq 0$, for $w \in I^{1,n}$, we obtain 
\[
\begin{aligned}
&\frac{P_n(\Phi;w,\varphi,\epsilon)}{P_n(\Phi;w,\psi,\epsilon)}\\
&=\frac{\sup\{\sum_{x\in E}e^{(S_{w,n}\varphi)(x)}|E~is~an\mathrm{~}(n,w,\epsilon;\Phi)\text{ separated subset of }X\}}{\sup\{\sum_{x\in E}e^{(S_{w,n}\psi)(x)}|E~is~an\mathrm{~}(n,w,\epsilon;\Phi)\text{ separated subset of }X\}}\\
&\leq\sup\left\{\frac{\sum_{x\in E}e^{(S_{w,n}\varphi)(x)}}{\sum_{x\in E}e^{(S_{w,n}\psi)(x)}} \Bigg| E~is~an~(n,w,\epsilon;\Phi)\text{ separated subset of }X\right\}\\
&\leq\sup\left\{\max_{x\in E}\frac{e^{(S_{w,n}\varphi)(x)}}{e^{(S_{w,n}\psi)(x)}}\Bigg|E~is~an~(n,w,\epsilon;\Phi)\text{ separated subset of }X\right\}\\
&\leq e^{(n+1)\|\varphi-\psi\|}.
\end{aligned}
\]
Hence
\[
P_n(\Phi;w,\varphi,\epsilon)\leq e^{(n+1)\|\varphi-\psi\|}P_n(\Phi;w,\psi,\epsilon).
\]
Moreover,
\[
P_n(\Phi;\varphi,\epsilon)\leq e^{(n+1)\|\varphi-\psi\|}P_n(\Phi;\psi,\epsilon).
\]
This implies
\[
P(\Phi;\varphi,\epsilon)-P(\Phi;\psi,\epsilon)\leq\|\varphi-\psi\|.
\]
Similarly, we can show that
\[
P(\Phi;\psi,\epsilon)-P(\Phi;\varphi,\epsilon)\leq\|\varphi-\psi\|.
\]
Therefore,
\[
|P(\Phi;\varphi,\epsilon)-P(\Phi;\psi,\epsilon)|\leq\|\varphi-\psi\|.
\]
If $P(\Phi,\cdot)$ is finite, taking the limit as $\epsilon \to 0$ yields
\[
|P(\Phi,\varphi)-P(\Phi,\psi)|\leq\|\varphi-\psi\|.
\]

$ (iv) $ For $w\in I^{1,n}$, by H\"{o}lder inequality, if $p\in[0,1]$ and $E$ is a finite subset of $X$,
we have
\[
\sum_{x\in E}e^{p(S_{w,n}\varphi)(x)+(1-p)(S_{w,n}\psi)(x)}\leq\left(\sum_{x\in E}e^{(S_{w,n}\varphi)(x)}\right)^p\left(\sum_{x\in E}e^{(S_{w,n}\psi)(x)}\right)^{1-p}.
\]
Thus
\[
P_{n}(\Phi;w,p\varphi+(1-p)\psi,\epsilon)\\\leq(P_{n}(\Phi;w,\varphi,\epsilon))^p(P_{n}(\Phi;w,\psi,\epsilon))^{1-p}.
\]
By H\"{o}lder inequality again
\begin{align*}
&\sum_{w \in I^{1,n}} P_{n}\left(\Phi;w, p\varphi+(1 - p)\psi,\epsilon\right)\\
\leq&\sum_{w \in I^{1,n}}\left(P_{n}\left(\Phi;w,\varphi,\epsilon\right)\right)^{p}\left(P_{n}\left(\Phi;w,\psi,\epsilon\right)\right)^{1 - p}\\
\leq&\left(\sum_{w \in I^{1,n}} P_{n}\left(\Phi;w,\varphi,\epsilon\right)\right)^{p}\left(\sum_{w \in I^{1,n}} P_{n}\left(\Phi;w,\psi,\epsilon\right)\right)^{1 - p}.
\end{align*}
Hence
\begin{align*}
&\frac{1}{\#\left(I^{1, n}\right)} \sum_{w \in I^{1,n}} P_{n}\left(\Phi;w, p\varphi+(1 - p)\psi,\epsilon\right)\\
\leq &\left(\frac{1}{\#\left(I^{1, n}\right)} \sum_{w \in I^{1,n}} P_{n}\left(\Phi;w,\varphi,\epsilon\right)\right)^{p}\left(\frac{1}{\#\left(I^{1, n}\right)} \sum_{w \in I^{1,n}} P_{n}\left(\Phi;w,\psi,\epsilon\right)\right)^{1 - p}.
\end{align*}
Thus
\[
P_{n}\left(\Phi;p\varphi+(1 - p)\psi,\epsilon\right)
\leq \left(P_{n}\left(\Phi;\varphi,\epsilon\right)\right)^{p}\left(P_{n}\left(\Phi;\psi,\epsilon\right)\right)^{1 - p}.
\]
Then we obtain
\[
P\left(\Phi; p\varphi+(1 - p)\psi,\epsilon\right)
\leq pP\left(\Phi;\varphi,\epsilon\right)+(1 - p)P\left(\Phi;\psi,\epsilon\right).
\]
If $P(\Phi,\cdot)$ is finite, taking the limit as $\epsilon \to 0$ yields 
\[
P\left(\Phi,p\varphi+(1 - p)\psi\right)\leq pP\left(\Phi,\varphi\right)+(1 - p)P\left(\Phi,\psi\right).
\]

$(v)$ is clear from the definition of topological pressure.

$(vi)$ For $w\in I^{1,n}$, note that $P_{n}\left(\Phi;w,\varphi+\psi,\epsilon\right)\leq P_{n}\left(\Phi;w,\varphi,\epsilon\right)P_{n}\left(\Phi;w,\psi,\epsilon\right).$ Thus
\begin{align*}
&\frac{1}{\#\left(I^{1, n}\right)}\sum_{w \in I^{1,n}}P_{n}\left(\Phi;w,\varphi+\psi,\epsilon\right)\\
\leq&\frac{1}{\#\left(I^{1, n}\right)}\sum_{w \in I^{1,n}}P_{n}\left(\Phi;w,\varphi,\epsilon\right)P_{n}\left(\Phi;w,\psi,\epsilon\right)\\
\leq& \#\left(I^{1, n}\right)\left(\frac{1}{\#\left(I^{1, n}\right)}\sum_{w \in I^{1,n}}P_{n}\left(\Phi;\varphi,\epsilon\right)\right)\left(\frac{1}{\#\left(I^{1, n}\right)}\sum_{w \in I^{1,n}}P_{n}\left(\Phi;\psi,\epsilon\right)\right).
\end{align*}
Then $P_{n}\left(\Phi;\varphi+\psi,\epsilon\right)\leq \#\left(I^{1, n}\right)P_{n}\left(\Phi;\varphi,\epsilon\right)P_{n}\left(\Phi;\psi,\epsilon\right).$
Hence 
\[
P\left(\Phi,\varphi+\psi\right)\leq P\left(\Phi,\varphi\right)+P\left(\Phi,\psi\right)+\limsup_{n \to \infty}\frac{\sum_{i=1}^{n}\log \#I^{(i)}}{n}.
\]

$(vii)$ If $a_1,\ldots,a_k$ are positive numbers with $\sum_{i= 1}^{k}a_i= 1 $ then $\sum_{i= 1}^{k}a_i^{c}\leq1$ if $c\geq1$, and $\sum_{i= 1}^{k}a_i^{c}\geq1$ if $c\leq1$. Therefore, for $w\in I^{1,n}$, if $E$ is an $(n,w,\epsilon;\Phi)$-separated subset of $X$ we have
\[
\sum_{x\in E}e^{c(S_{w,n}\varphi)(x)}\leq\left(\sum_{x\in E}e^{(S_{w,n}\varphi)(x)}\right)^{c}\text{ if }c\geq1
\]
and
\[
\sum_{x\in E}e^{c(S_{w,n}\varphi)(x)}\geq\left(\sum_{x\in E}e^{(S_{w,n}\varphi)(x)}\right)^{c}\text{ if }c\leq1.
\]
Therefore $P_{n}\left(\Phi;w,c\varphi,\epsilon\right)\leq\left(P_{n}\left(\Phi;w,\varphi,\epsilon\right)\right)^{c}$ if $c\geq1$ and $P_{n}\left(\Phi;w,c\varphi,\epsilon\right)\geq\left(P_{n}\left(\Phi;w,\varphi,\epsilon\right)\right)^{c}$ if $c\leq1$.

\textbf{Case 1}, if $c\geq1$, then
\begin{align*}
&\frac{1}{\#(I^{1,n})}\sum_{w \in I^{1,n}}P_{n}\left(\Phi;w,c\varphi,\epsilon\right)\\
\leq&\frac{1}{\#(I^{1,n})}\sum_{w \in I^{1,n}}\left(P_{n}\left(\Phi;w,\varphi,\epsilon\right)\right)^{c}\\
\leq&\frac{1}{\#(I^{1,n})}\left(\sum_{w \in I^{1,n}}P_{n}\left(\Phi;w,\varphi,\epsilon\right)\right)^{c}\\
=&\left(\frac{1}{\#(I^{1,n})}\right)^{1 - c}\left(\frac{1}{\#(I^{1,n})}\sum_{w \in I^{1,n}}\left(P_{n}\left(\Phi;w,\varphi,\epsilon\right)\right)\right)^{c}.
\end{align*}

\textbf{Case 2}, if $c\leq1$, then
\begin{align*}
&\frac{1}{\#(I^{1,n})}\sum_{w \in I^{1,n}}P_{n}\left(\Phi;w,c\varphi,\epsilon\right)\\
\geq &\frac{1}{\#(I^{1,n})}\sum_{w \in I^{1,n}}\left(P_{n}\left(\Phi;w,\varphi,\epsilon\right)\right)^{c}\\
\geq & \frac{1}{\#(I^{1,n})}\left(\sum_{w \in I^{1,n}}P_{n}\left(\Phi;w,\varphi,\epsilon\right)\right)^{c}\\
=&\left(\frac{1}{\#(I^{1,n})}\right)^{1 - c}\left(\frac{1}{\#(I^{1,n})}\sum_{w \in I^{1,n}}\left(P_{n}\left(\Phi;w,\varphi,\epsilon\right)\right)\right)^{c}
\end{align*}
which implies the desired result.

$(viii)$ Since $-\vert\varphi\vert\leq\varphi\leq\vert\varphi\vert$, by $(i)$, $P(\Phi,-\vert\varphi\vert)\leq P(\Phi,\varphi)\leq P(\Phi,\vert\varphi\vert)$. From $(vii)$ we have 
\[
-P\left(\Phi,\vert\varphi\vert\right)-2\liminf_{n \to \infty}\frac{\sum_{i=1}^{n}\log \#I^{(i)}}{n}\leq P\left(\Phi,-\vert\varphi\vert\right).
\] 
Thus, 
\[
-2\liminf_{n \to \infty}\frac{\sum_{i=1}^{n}\log \#I^{(i)}}{n}-P\left(\Phi,\vert\varphi\vert\right)\leq P\left(\Phi,\varphi\right)\leq P\left(\Phi,\vert\varphi\vert\right).
\]
\end{proof}

\end{theorem}

An NAIFS $(X, \Phi)$ of continuous maps on a compact metric space $(X, d)$ is said to be equi-continuous, if for every $\epsilon>0$ there exists $\delta>0$ such that the implication $d(x, y)<\delta \Rightarrow d\left(f_{i}^{(j)}(x), f_{i}^{(j)}(y)\right)<\epsilon$ holds for every $x, y \in X, j \geq 1$ and $i \in I^{(j)}$.

\begin{theorem}
Let \((X,\Phi)\) be an equicontinuous NAIFS of a compact metric space \((X,d)\) and \(\psi = c\), where \(c\) is a constant. Then, \(P(\Phi^n,n\psi)=nP(\Phi,\psi)\) for all \(n\geq 1\).
\begin{proof}
We know that \(h(\Phi^n)=nh(\Phi)\) for all \(n\geq 1\) (see \cite{ghane2019}). Since \(\psi = c\) and according to \((v)\) of Theorem \ref{th1},
\[
P(\Phi^n,n\psi)=nc + h(\Phi^n)=nc + nh(\Phi)=n(c + h(\Phi))=nP(\Phi,\psi).
\]
\end{proof}
\end{theorem}

\begin{theorem}\label{pow1}
Let $ (X,\Phi)$ be an NAIFS of a compact metric space $(X,d)$ and let $\psi:X\to \mathbb{R}^{+}$ be a non-negative and continuous map. Then for every $n \geq 1$,
\[
P\left(\Phi^{n},\psi\right) \leq nP\left(\Phi,\psi\right).
\]
\begin{proof}
For any $w=w_{1}w_{2} \cdots w_{mn} \in I^{1, mn}$, denote $w^{*}=w_{1}^{*} w_{2}^{*} \cdots w_{m}^{*} \in I^{1,m}_{*},$ where $w_{i}^{*}=w_{(i-1)n+1} w_{(i-1)n+2} \cdots w_{in}$ and $1 \leq i \leq m.$ Then, for any compact subset $K \subset
 X$, let $F \subset X$ be a $(mn,w,\epsilon;\Phi)$-spanning set for $K$ of cardinality $r_{mn}(K;w,\epsilon,\Phi)$. Therefore, for each $x \in K$, there exits $y \in F$ such that
\[
d_{w,mn}(x,y)=\max_{0 \leq i \leq mn}d(f_{w}^{i}(x),f_{w}^{i}(y)) \leq \epsilon.
\]
It is obvious that $d_{w}(x,y) \geq d_{w*}(x,y),$ then $F$ is also a $(m,w^{*},\epsilon;\Phi^{n})$-spanning set for $K$. So, $r_{m}(K;w^{*},\epsilon,\Phi^{n}) \leq r_{mn}(K;w,\epsilon,\Phi)$.  And then we observe that
\[
S_{w*,n}\psi(x)=\sum_{i=0}^{n}\psi\left(f_{w^{*}}^{1,i}(x)\right)\leq\sum_{i=0}^{mn}\psi\left(f_{w}^{1,i}(x)\right)=S_{w,mn}\psi\left(x\right),
\]
since $\psi:X\to\mathbb{R}^{+}$ is a non-negative map. Then
\[
Q_{n}(\Phi^{n};w^{*},\psi,\epsilon)\leq \sum_{x\in F}e^{S_{w^{*},n}\psi(x)} \leq \sum_{x\in F}e^{S_{w,mn}\psi(x)}.
\]
Hence,
\[
Q_{n}(\Phi^{n};w^{*},\psi,\epsilon)\leq Q_{mn}(\Phi;w,\psi,\epsilon).
\]
It follows that
\[
\begin{aligned}
& P\left(\Phi,\psi\right) \\
=& \lim_{\epsilon \rightarrow 0} \limsup_{n \rightarrow \infty} \frac{1}{n}\log \frac{1}{\#\left(I^{1, n}\right)} \sum_{w \in I^{1, n}} Q_{n}(\Phi;w,\psi,\epsilon) \\
\geq & \lim_{\epsilon \rightarrow 0} \limsup_{m \rightarrow \infty} \frac{1}{mn}\log \frac{1}{\#\left(I^{1, mn}\right)} \sum_{w \in I^{1, mn}} Q_{mn}(\Phi;w,\psi,\epsilon) \\
\geq & \lim_{\epsilon \rightarrow 0} \limsup_{m \rightarrow \infty} \frac{1}{mn} \log \frac{1}{\#\left(I^{1, m}_{*} \right)} \sum_{w^{*} \in I^{1,m}_{*}}Q_{n}(\Phi^{n};w^{*},\psi,\epsilon) \\
=& \frac{1}{n} \lim_{\epsilon \rightarrow 0} \limsup_{m \rightarrow \infty} \frac{1}{m} \log \frac{1}{\#\left(I^{1, m}_{*} \right)} \sum_{w^{*} \in I^{1,m}_{*}} Q_{n}(\Phi^{n};w^{*},\psi,\epsilon)= \frac{1}{n} P\left(\Phi^{n},\psi\right).
\end{aligned}
\]
Hence, $ P\left(\Phi^{n},\psi\right) \leq n \cdot  P\left(\Phi,\psi\right)$ which completes the proof.
\end{proof}
\end{theorem}

Let us take an NAIFS $(X, \Phi)$ in which $X$ is a compact metric space and $\Phi$ consists of a sequence $\left\{\Phi^{(j)}\right\}_{j \geq 1}$ of collections of maps, where $\Phi^{(j)}=\left\{f_i^{(j)}\right.$ : $X \rightarrow X\}_{i \in I^{(j)}}$ and $I^{(j)}$ is a non-empty finite index set for all $j \geq 1$. For each $k \geq 1$ we will denote by $\left(X, \Phi_k\right)$ the NAIFS composed of the sequence $\left\{\Phi^{(j)}\right\}_{j \geq k}$.

\begin{lemma}\label{mono}
Let $ (X,\Phi)$ be an NAIFS on a compact metric space $(X,d)$. Thus, for any $ 1 \leq i <j < \infty $ and every $\psi \in C(X,\mathbb{R})$, we have $P\left(\Phi_{i},\psi\right) \leq P\left(\Phi_{j},\psi\right)$.
\begin{proof}
This theorem can be proved by demonstrating that $P\left(\Phi_{i},\psi\right) \leq P\left(\Phi_{i+1},\psi\right)$ for all $i \geq 1$.
For any $ \epsilon >0 $, due to the compactness of $X$, we can choose a finite subset $ F' \subset X$ such that 
\[
X  \subset \bigcup_{y' \in F'}B_{d}(y',\epsilon).
\]
For $w=w_i w_{i+1} \cdots w_{i+n-1} \in I^{i, n}$, put $w^{\prime}=w_{i+1} \cdots w_{i+n-1} \in I^{i+1, n-1}$. Let $ F $ be an $ (n-1,w',\epsilon;\Phi_{i+1}) $-spanning set of $ X $, i.e.,
\[
X \subset \bigcup_{y \in F}\bigcap_{j=0}^{n-1}f_{w'}^{i+1,-j}(B_{d}(f_{w'}^{i+1,j}(y),\epsilon)).
\]
Then
\[
\begin{aligned}
X &\subset \left(\bigcup_{y' \in F'}B_{d}(y',\epsilon)\right) \bigcap \left(f_{w_i}^{i,-1}\left(\bigcup_{y \in F}\bigcap_{j=0}^{n-1}f_{w'}^{i+1,-j}(B_{d}(f_{w'}^{i+1,j}(y),\epsilon))\right)\right)\\
& =\left(\bigcup_{y' \in F'}B_{d}(y',\epsilon)\right) \bigcap \left(\bigcup_{y \in F}\bigcap_{j=1}^{n}f_{w}^{i,-j}(B_{d}(f_{w'}^{i+1,j-1}(y),\epsilon))\right)\\
&=\bigcup_{y' \in F'}\bigcup_{y \in F}\left(B_{d}(y',\epsilon) \bigcap \bigcap_{j=1}^{n}f_{w}^{i,-j}(B_{d}(f_{w'}^{i+1,j-1}(y),\epsilon)) \right).
\end{aligned}
\]
Without loss of  generality, we assume that
\[
B_{d}(y',\epsilon) \bigcap \bigcap_{j=1}^{n}f_{w}^{i,-j}(B_{d}(f_{w'}^{i+1,j-1}(y),\epsilon)) \neq \emptyset.
\]
One can select, for every $y' \in  F'$ and $ y \in F $, exactly one point $ x^* \in B_{d}(y',\epsilon) \bigcap \bigcap_{j=1}^{n}f_{w}^{i,-j}(B_{d}(f_{w'}^{i+1,j-1}(y),\epsilon))$ to construct a set $ F^{*} $. For any $ x \in B_{d}(y',\epsilon) $, we have  
\[
d(x^*,x) \leq d(x^*,y')+d(x,y') < 2\epsilon.
\]
Hence, $ x \in  B_{d}(x^*,2\epsilon)$ which follows $B_{d}(y',\epsilon) \subset  B_{d}(x^*,2\epsilon) $.
Take $ x \in f_{w}^{i,-j}(B_{d}(f_{w'}^{i+1,j-1}(y),\epsilon)) $ for any $ j \in \{1,\ldots,n \} $, we have 
$ d(f_{w}^{i,j}(x),f_{w'}^{i+1,j-1}(y)) < \epsilon $. Since $ x^* \in f_{w}^{i,-j}(B_{d}(f_{w'}^{i+1,j-1}(y),\epsilon))$, then 
\[
d(f_{w}^{i,j}(x^*),f_{w}^{i,j}(x)) \leq d(f_{w}^{i,j}(x^*),f_{w'}^{i+1,j-1}(y))+d(f_{w}^{i,j}(x),f_{w'}^{i+1,j-1}(y)) < 2\epsilon.
\]
Hence, $ x \in  f_{w}^{i,-j}B_{d}(f_{w}^{i,j}(x^*),2\epsilon)$ which follows
\[
f_{w}^{i,-j}(B_{d}(f_{w'}^{i+1,j-1}(y),\epsilon))  \subset  f_{w}^{i,-j}(B_{d}(f_{w}^{i,j}(x^*),2\epsilon)).
\]
Therefore
\[
B_{d}(y',\epsilon) \bigcap \bigcap_{j=1}^{n}f_{w}^{i,-j}(B_{d}(f_{w'}^{i+1,j-1}(y),\epsilon))  \subset \bigcap_{j=0}^{n}f_{w}^{i,-j}B_{d}(f_{w}^{i,j}(x^*),2\epsilon).
\]
It is obvious that $ F^* $ is an $ (n,w,2\epsilon;\Phi_{i}) $-spanning set of $ X$, then we have
\[
r_{n}\left(w, 2\epsilon, \Phi_{i}\right) \leq \#(F') \cdot r_{n-1}\left(w', \epsilon, \Phi_{i+1}\right).
\]
If $\delta=\sup\{|\psi(x)-\psi(y)||d(x,y) \leq \epsilon\}$, then we have 
\[
\begin{aligned}
&Q_{n}(\Phi_{i};w,\psi,2\epsilon )\\
&\leq\sum_{x^{*}\in F^{*}}\exp\left(\sum_{j=1}^{n}\psi(f_{w}^{i,j}(x^{*}))+\psi(x^{*})\right)\\
&\leq e^{\Vert \psi \Vert}\sum_{x^{*}\in F^{*}}\exp\left(\sum_{j=1}^{n}\psi(f_{w}^{i,j}(x^{*}))\right)\\
&=e^{\Vert \psi \Vert}\sum_{x^{*}\in F^{*}}\exp\left(\sum_{j=1}^{n}\left(\psi(f_{w}^{i,j}(x^{*}))-\psi(f_{w'}^{i+1,j-1}(y))+\psi(f_{w'}^{i+1,j-1}(y))\right)\right)\\
&=e^{\Vert \psi \Vert}\sum_{y'\in F'}\sum_{y\in F}\exp\left(\sum_{j=1}^{n}\psi(f_{w'}^{i+1,j-1}(y))+\sum_{j=1}^{n}\left(\psi(f_{w}^{i,j}(x^{*}))-\psi(f_{w'}^{i+1,j-1}(y))\right)\right)\\
&\leq e^{\Vert \psi \Vert}\sum_{y'\in F'}\sum_{y\in F}\exp\left(\sum_{j=1}^{n}\psi(f_{w'}^{i+1,j-1}(y))+n\delta\right)\\
&=e^{\Vert \psi \Vert+n\delta}\#(F')\sum_{y\in F}\exp\left(\sum_{j=1}^{n}\psi(f_{w'}^{i+1,j-1}(y))\right).
\end{aligned}
\]
Hence,
\[
Q_{n}(\Phi_{i};w,\psi,2\epsilon ) \leq  e^{\Vert \psi \Vert+n\delta}\#(F')Q_{n-1}(\Phi_{i+1};w',\psi,\epsilon ).
\]
It follows that  
\[
\begin{aligned}
P\left(\Phi_{i};\psi,2\epsilon \right)&=\limsup_{n \to \infty }\frac{1}{n} \log Q_n(\Phi_{i};\psi,2\epsilon)\\
& =\limsup_{n \to \infty }\frac{1}{n}\log\frac{1}{\#(I^{i,n})}\sum_{w\in I^{i,n}}Q_n(\Phi_{i};w,\psi,2\epsilon) \\
& \leq \limsup_{n \to \infty }\frac{1}{n}\log\frac{1}{\#(I^{i,n})}\sum_{w\in I^{i,n}}e^{\Vert \psi \Vert+n\delta}\#(F')Q_{n-1}(\Phi_{i+1};w',\psi,\epsilon ) \\
& = \limsup_{n \to \infty }\frac{1}{n}\log\frac{\#(I^{(i)})}{\#(I^{i,n})}\sum_{w'\in I^{i+1,n-1}}Q_{n-1}(\Phi_{i+1};w',\psi,\epsilon )+\delta\\
& = \limsup_{n \to \infty }\frac{1}{n}\log\frac{1}{\#(I^{i+1,n-1})}\sum_{w'\in I^{i+1,n-1}}Q_{n-1}(\Phi_{i+1};w',\psi,\epsilon)+\delta\\
& = P(\Phi_{i+1};\psi,\epsilon)+\delta.
\end{aligned}
\]
Since $\delta\to0$ as $\epsilon\to0$, now by letting $ \epsilon \rightarrow 0 $, we get $P\left(\Phi_{i},\psi\right) \leq P\left(\Phi_{i+1},\psi\right) $ which completes the proof.
\end{proof}
\end{lemma}

\section{Main results}\label{sec5}

\subsection{Topological conjugacy}

Let $\left(X, \Phi\right)$ and $\left(Y, \Psi\right)$ be NAIFSs of continuous self-maps on the compact metric spaces $(X,d)$ and $(Y,\rho)$, respectively. Assume that 
$
\Phi^{(j)} = \left\{ f^{(j)}_i : X \rightarrow X \right\}
$
and 
$
\Psi^{(j)} = \left\{ g^{(j)}_i : Y \rightarrow Y \right\}
$
, with $i \in I^{(j)}$ and $j \in \mathbb{N}$ (i.e., the index sets for $\Phi$ and $\Psi$ are identical). Let $\pi$ be a surjective and continuous map from $X$ into $Y$ such that the diagram commutes 
$$
\xymatrix { 
X \ar[rrr]^{f_{i}^{(j)}} \ar[d]_{\pi} & & & X \ar[d]^{\pi}\\
Y \ar[rrr]_{g_{i}^{(j)}} & & & Y
}
$$
(i.e., $\pi \circ f_{i}^{(j)}=g_{i}^{(j)} \circ \pi$ for every $i \in I^{(j)}$). $(Y, \Psi)$ is (topological) semiconjugate to $(X, \Phi)$, with $ \pi $ as a factor map. If $\pi$ is a homeomorphism, then $(X, \Phi)$ is (topological) conjugate to $(Y, \Psi)$ and $ \pi $ is a (topological) conjugacy.

The proofs of Theorem \ref{con1} and Corollary \ref{con2} closely follow Theorem 2.1 and Corollary 2.1 in \cite{Zhao}, so the details are omitted here. Interested readers may refer to \cite{Zhao} for a detailed exposition of the argument.
\begin{theorem}\label{con1}
Let $\left(X, \Phi\right)$ and $\left(Y, \Psi\right)$ be NAIFSs of continuous self-maps on the compact metric spaces $(X,d)$ and $(Y,\rho)$, respectively. Assume that $\Phi^{(j)}=\left\{f^{(j)}_{i}: X \rightarrow X \right\}$, $\Psi^{(j)}=\left\{g^{(j)}_{i}: Y \rightarrow Y\right\}$, with $i \in I^{(j)}, j \in \mathbb{N}$ and that $(X,\Phi)$ is semiconjugate to $(Y,\Psi)$. Let \(\varphi: Y\rightarrow\mathbb{R}\) be a continuous map, and the factor map be \(\pi: X\rightarrow Y\), then
\[
P(\Psi,\varphi) \leq P(\Phi,\varphi\circ\pi).
\]
\end{theorem}

\begin{corollary}\label{con2}
If the system $(X,\Phi)$ is conjugate to $(Y,\Psi),$ then
\[
P(\Psi,\varphi) = P(\Phi,\varphi\circ\pi).
\]
\end{corollary}

To prove Theorem \ref{con}, we will primarily follow the approach used in \cite{Bowen1} and \cite{KS}. First, we introduce the concept of the topological sup-entropy of $\Phi$, which will be useful in the subsequent analysis.

Let \((X,\Phi)\) be an equicontinuous NAIFS of a compact metric space \((X,d)\) and \(Y\) be a nonempty subset of \(X\). For each \(n\geq1\) the function \(d^{*}_{n}\) given by 
\[
d^{*}_{n}(x,y)=\sup_{i\in \mathbb{N}}\max_{w \in I^{i,n}}\max_{0\leq t \leq  n}d\left(f_{w}^{i,t}(x),f_{w}^{i,t}(y)\right)
\]
is a distance on \(X\). Since $\Phi$ is assumed to be composed of equicontinuous maps, $d_n^{*}$ is equivalent to $d$. Hence,  the metric space $(X, d_n^{*})$ is compact. A subset \(E^{*}\) of the space \(X\) is called \((n,\epsilon)^{*}\)-separated if for any two distinct points \(x,y\in E\), \(d^{*}_{n}(x,y)>\epsilon\). A set \(F^{*}\subset X\) \((n,\epsilon)^{*}\)-spans another set \(K\subset X\) provided that for each \(x\in K\) there is \(y\in F^{*}\) for which \(d^{*}_{n}(x,y)\leq\epsilon\).

We define \(s_n^{*}(\Phi;Y;\epsilon)\) as the maximal cardinality of an \((n,\epsilon)^{*}\)-separated set in \(Y\) and \(r_n^{*}(\Phi;Y;\epsilon)\) as the minimal cardinality of a set in \(Y\) which \((n,\epsilon)^{*}\)-spans \(Y\). 

Exactly as in the proof of Lemma 3.1 in \cite{ghane2019} one can show that 
\[
r_{n}^{*}(\Phi;Y;\epsilon)\leq s_{n}^{*}(\Phi;Y;\epsilon)\leq r_{n}^{*}(\Phi;Y;\epsilon/2)
\]
and so we can define 
\[
H(\Phi;Y)=\lim_{\varepsilon\to0}\limsup_{n\to\infty}\frac{1}{n}\log s_{n}^{*}(\Phi;Y;\varepsilon)=\lim_{\varepsilon\to0}\limsup_{n\to\infty}\frac{1}{n}\log r_{n}^{*}(\Phi;Y;\varepsilon).
\]
We also put \(H(\Phi;\varnothing)=0\).
The quantity \(H(\Phi;Y)\) is said to be the topological sup-entropy of \(\Phi\) on the set \(Y\).

\begin{theorem}\label{con}
Let $\left(X, \Phi\right)$ and $\left(Y, \Psi\right)$ be equicontinuous NAIFSs on the compact metric spaces $(X,d)$ and $(Y,\rho)$, respectively. We assume that $\Phi^{(j)}=\left\{f^{(j)}_{i}: X \rightarrow X \right\}$, $\Psi^{(j)}=\left\{g^{(j)}_{i}: Y \rightarrow Y\right\}, i \in I^{(j)}, j \in \mathbb{N}$ and $(X,\Phi)$ is semiconjugate with $(Y,\Psi)$. Let \(\varphi: Y\rightarrow\mathbb{R}\) be a continuous map, and the factor map be \(\pi: X\rightarrow Y\), then 
\[
P(\Phi,\varphi\circ\pi)\leq P(\Psi,\varphi)+\sup_{y\in Y}H(\Phi;\pi^{-1}(y)).
\] 
\begin{proof}
In the following,  for any $\epsilon> 0$, we set
\[
\mathrm{var}(\varphi\circ\pi,\epsilon)=\sup_{x,y\in X}\{|\varphi\circ\pi(x)-\varphi\circ\pi(y)|:d(x,y)<\epsilon\}.
\]
We can assume that $a= \sup_{y \in Y}H(\Phi; \pi ^{-1}(y) ) < \infty$, since if $a=\infty$, there is nothing to prove. 

Let $r_{n}^{*}(\Phi;\pi^{-1}(y);\epsilon)$ be the minimum cardinality of a subset of $\pi^{-1}(y)$ which $(n,\epsilon)^{*}$-spans $\pi^{-1}(y)$ with respect to $\Phi$. As $r_n^*(\Phi;\pi^{-1}(y);\epsilon)$ decreases with $\epsilon>0$, we have
\[
H(\Phi;\pi ^{- 1}( y) ) \geq H(\Phi;\pi ^{-1}(y),\epsilon),
\]
where
\[
H(\Phi;\pi^{-1}(y),\epsilon)=\limsup_{n\to\infty}\frac{1}{n}r_{n}^{*}(\Phi;\pi^{-1}(y);\epsilon).
\]

Fix $\epsilon>0$ and choose $\alpha>0$, for any $y\in Y$, we can choose an integer $m(y)$ such that
\begin{equation}\label{fa1}
a+\alpha\geq H(\Phi;\pi^{-1}(y),\epsilon)+\alpha \geq\frac1{m(y)}\log r_{m(y)}^{*}(\Phi;\pi^{-1}(y);\epsilon).
\end{equation}
In order to prove the theorem, we choose a set $F_y^*\subset\pi^{-1}(y)$ with the smallest possible cardinality so that it $(m(y),\epsilon)^*$-spans $\pi^{-1}(y)$ with respect to $\Phi$. From \eqref{fa1} we have
\begin{equation}\label{fa2}
a+\alpha\geq\frac{1}{m(y)}\log\#(F_{y}^{*}).
\end{equation}
Let $D^{*}_{m}(z, 2 \epsilon, \Phi)$ be the neighbourhood of $z$ defined by
\[
D^{*}_{m}(z, 2 \epsilon, \Phi):=\left\{c \in X: d^{*}_{m}(c, z)<2 \epsilon\right\},
\]
where 
\[
d^{*}_{m}(c,z)=\sup_{i \in \mathbb{N}}\max_{w \in I^{i,m}}\max_{0\leq t \leq  m}d\left(f_{w}^{i,t}(c),f_{w}^{i,t}(z)\right)< 2\epsilon,
\]
and also denote
\[
U_{y}=\bigcup_{z \in F_y^{*}} D^{*}_{m(y)}(z, 2 \epsilon, \Phi).
\]
Then $U_{y}$ is an open neighborhood of $\pi^{-1}(y)$ and
\[
\left(X \backslash U_{y}\right) \cap \bigcap_{r>0} \pi^{-1}\left(\overline{B_{r}(y)}\right)=\emptyset,
\]
where $B_r (y)=\{y' \in Y: \rho(y', y)<r\}$. By the finite intersection property, there is a $W_y = B_{r_y} (y)$ such that $U_y \supset \pi^{-1}(W_y)$. Since $ Y $ is compact, there exists $\{y_1, y_2, \ldots, y_p\}$ such that  $\{W_{y_1}, W_{y_2}, \ldots, W_{y_p}\}$ covers $Y$. Let $\delta<\delta_{1}/2$ where $\delta_{1}$ is a Lebesgue number of this cover.

By the definition of $P(\Psi, \varphi)$, we have
$
P(\Psi, \varphi) \geq P(\Psi;\varphi,\delta),
$
where
\[
P(\Psi;\varphi,\delta)=\limsup_{n\to\infty}\frac{1}{n}\log P_{n}(\Psi;\varphi,\delta).
\]
Fix $ \beta>0 $, we can choose an integer $ N $ such that for any $ n>N $,
\begin{equation}\label{fa3}
\begin{aligned}
P(\Psi, \varphi)+\beta 
\geq P(\Psi;\varphi,\delta)+\beta 
> &\frac{1}{n} \log \frac{1}{\#(I^{1,n})} \sum_{w\in I^{1,n} } P_n(\Psi;w,\varphi,\delta)\\
\geq & \frac{1}{n} \log \frac{1}{\#(I^{1,n})} \sum_{w\in I^{1,n} } \sum_{y \in E_{w, n}} e^{S_{w, n} \varphi(y)},
\end{aligned}
\end{equation} 
where $ E_{w, n} $ is an $(n,w,\delta;\Psi)$-maximal separated set of $Y$.
Since $ Y $ is compact, we can assume the $\delta$ is small enough such that for any $ x,y \in Y $,
\begin{equation}\label{fa4}
\sup_{\rho(x,y) <\delta}|\varphi(x)-\varphi(y)|<\epsilon.
\end{equation}

Take $y\in Y$ and let $c_{0}(y)\in \{y_{1},...,y_{p}\}$ be such that $W_{c_{0}(y)}\supset\overline{B}_{\delta}(y)$ and define $t_{0}(y)=0$. Next let $t_{1}(y)=m(c_{0}(y))$ and let $c_{1}(y)\in\{y_{1},...,y_{p}\}$ be such that $W_{c_1(y)}\supset\overline{B}_\delta(g_{w}^{1,t_1(y)}(y)).$ Similarly, assuming that $t_0(y),...,t_k(y)$ and $c_0(y),...,c_k(y)$ are already defined we define $t_{k+1}(y)=t_k(y)+m(c_k(y))$ and let $c_{k+1}(y)\in\{y_1,...,y_p\}$ be such that $W_{c_{k+1}(y)}\supset\overline{B}_{\delta}(g_{w}^{1,t_{k+1}(y)}(y)).$ Finally, let $l=l(y)$ be such that
\begin{equation}\label{fa5}
\sum_{s=0}^{l-1}m(c_{s}(y))=t_{l}(y)<n\leq t_{l}(y)+m(c_{l}(y)).
\end{equation}

For each $y \in E_{w, n}$, $x_0 \in F^{*}_{c_{0}(y)}, \ldots, x_l \in F^{*}_{c_{l}(y)}$, we consider the following set
\begin{align*}
V_{w}\left(y ; x_{0}, \ldots, x_{l}\right)
=\big\{ & x \in X :d\left(f_{w}^{1,t+t_{s}(y)}(x), f_{w}^{t_{s}(y)+1,t}\left(x_{s}\right)\right)<2 \epsilon, \\
&  0 \leq t < m\left(c_{s}(y)\right), 0 \leq s \leq l~\text{and}~t_{l}(y)+t \leq n \big\}.
\end{align*}
Then we claim that
\begin{enumerate}[label=(\alph*)]
    \item \label{a} $\mathcal{V}_{w} = \{ V_{w}( y; x_0, \ldots , x_s) : y \in E_{w, n}, x_s\in F^{*}_{c_{s}(y)}, 0\leq s\leq l( y) \}$ is a cover of $X.$
    
    \item \label{b} any $(n,w,4\varepsilon;\Phi)$-separated set intersects each element of $\mathcal{V}_{w}$ in at most one point.
\end{enumerate}
In order to prove \ref{a}, let $x \in X$, since $E_{w, n}$ is also an $(n, w, \delta; \Psi)$-spanning set of $Y$, then there exists  $y \in E_{w, n}$ such that
\[
\rho_{w,n}(y, \pi(x)) = \max_{ 0 \leq i \leq n} \rho(g_{w}^{1,i}(y), g_{w}^{1,i}(\pi(x))) < \delta.
\]
For any $0 \leq s \leq l(y)$, 
\[
\pi \circ f_{w}^{1,t_{s}(y)}(x) = g_{w}^{1,t_{s}(y)}(\pi(x)) \in W_{c_{s}(y)}.
\]
This implies that there exists $x_{s} \in F^{*}_{c_{s}(y)}$ such that 
\[
d\left(f_{w}^{t_{s}(y)+1,t}(f_{w}^{1,t_{s}(y)}(x)), f_{w}^{t_{s}(y)+1,t}\left(x_{s}\right)\right)<2 \epsilon
\]
for all $0 \leq t < m(c_{s}(y)), 0 \leq s \leq l(y)$ and $t_{l}(y)+t \leq n$. Thus $x \in V_{w}(y; x_0, \ldots, x_l)$. This proves \ref{a}.

If $z, c \in V_{w}(y, x_0, \ldots, x_l)$, then 
\[
\begin{aligned}
&d\left(f_{w}^{1,t+t_{s}(y)}(z), f_{w}^{1,t+t_{s}(y)}(c)\right)\\
\leq & d\left(f_{w}^{1,t+t_{s}(y)}(z), f_{w}^{1+t_{s}(y),t}\left(x_{s}\right)\right) + d\left(f_{w}^{1,t+t_{s}(y)}(c), f_{w}^{1+t_{s}(y),t}\left(x_{s}\right)\right) < 4 \epsilon,
\end{aligned}
\]
for each $0 \leq t < m(c_{s}(y))$ and $0 \leq s \leq l(y)$. This proves \ref{b}.

For any $y \in E_{w, n}$, by the definition of $V_{w}(y; x_0, \ldots, x_l)$, we can denote 
\[
\mathcal{V}_{w,y}:=\left\{V_{w}\left(y; x_{0}, \ldots, x_{l}\right): x_{s} \in F_{c_{s}(y)}^{*} \text{ for all } 0 \leq s \leq l\right\}.
\]
Set $\mathcal{V}_{w}= \bigcup_{y \in E_{w, n}} \mathcal{V}_{w,y}$, where the cardinality of $\mathcal{V}_{w,y}$ is
\[
\#\left(\mathcal{V}_{w,y}\right)=\prod_{s=0}^{l=l(y)} \#\left(F_{c_{s}(y)}^{*}\right).
\]
By \eqref{fa2} and \eqref{fa5}, 
\begin{equation}\label{fa6}
\begin{aligned}
\#(\mathcal{V}_{w,y})
&=\exp\left(\sum_{s=0}^{l(y)}\log \#\left(F_{c_{s}(y)}^{*}\right)\right)\\
&\leq \exp\left((a+\alpha)\sum_{s=0}^{l(y)}m(c_{s}(y))\right)\\
&=\exp\left((a+\alpha)\left(\sum_{s=0}^{l(y)-1}m(c_{s}(y))+m(c_{l}(y))\right)\right)\\
&\leq \exp((a+\alpha)(n+M)),
\end{aligned}
\end{equation}
where $M=\max\{m(y_1),\ldots,m(y_p)\}.$
For any $(n, w, 4\epsilon;\Phi)$-separated set $H_{w, n}$ of $X$, we want to estimate the upper bound of
\[
\sum_{x \in H_{w, n}} e^{S_{w, n} \varphi \circ \pi(x)}.
\]
By \ref{a}, for each $x \in H_{w, n}$, there exist $y \in E_{w, n}$ and $x_0 \in F_{c_{0}(y)}^{*}, \ldots, x_l \in F_{c_{l}(y)}^{*}$ such that $x \in V_{w}(y; x_0, \ldots, x_l)$.
Furthermore,
\begin{align*}
& S_{w, n} \varphi \circ \pi(x) \\
=& \sum_{t=0}^{n} \varphi\left(\pi \circ f_{w}^{1,t}(x)\right) \\
=& \sum_{s=0}^{l-1} \sum_{t=0}^{m\left(c_{s}(y)\right) - 1} \varphi\left(\pi \circ f_{w}^{1,t_s(y)+t}(x)\right)+ \sum_{t=0}^{n-t_l(y)} \varphi\left(\pi \circ f_{w}^{1,t_l(y)+t}(x)\right)\\
=& \sum_{s=0}^{l-1}\sum_{t=0}^{m\left(c_{s}(y)\right) - 1}\varphi\left(\pi \circ f_{w}^{t_s(y)+1,t}\left(x_{s}\right)\right)+ \sum_{t=0}^{n-t_l(y)} \varphi\left(\pi \circ f_{w}^{t_l(y)+1,t}(x_l)\right)+\\
&\sum_{s=0}^{l-1}\sum_{t=0}^{m\left(c_{s}(y)\right) - 1}\left(\varphi\left(\pi \circ f_{w}^{1,t_s(y)+t}(x)\right)-\varphi\left(\pi \circ f_{w}^{t_s(y)+1,t}\left(x_{s}\right)\right)\right)+\\
&\sum_{t=0}^{n-t_l(y)} \left(\varphi\left(\pi \circ f_{w}^{1,t_l(y)+t}(x)\right)- \varphi\left(\pi \circ f_{w}^{t_l(y)+1,t}(x_l)\right)\right).
\end{align*}
Since $d\left(f_{w}^{1,t_{s}(y)+t}(x), f_{w}^{t_{s}(y)+1,t}\left(x_{s}\right)\right)<2 \epsilon$ for any \(0 \leq t<m(c_{s}(y)\), \(0 \leq s \leq l\) and \(t_{l}(y)+ t \leq n\), then 
\[
\varphi\left(\pi \circ f_{w}^{1,t_s(y)+t}(x)\right)-\varphi\left(\pi \circ f_{w}^{t_s(y)+1,t}\left(x_{s}\right)\right) \leq \operatorname{var}(\varphi \circ \pi, 2 \epsilon).
\]
Therefore,
\begin{align*}
S_{w, n} \varphi \circ \pi(x) &\leq \sum_{s=0}^{l-1}\sum_{t=0}^{m\left(c_{s}(y)\right)-1}\varphi\left(\pi \circ f_{w}^{t_s(y)+1,t}\left(x_{s}\right)\right)+\sum_{t=0}^{n-t_l(y)} \varphi\left(\pi \circ f_{w}^{t_l(y)+1,t}(x_l)\right)\\
&+(n+1) \operatorname{var}(\varphi \circ \pi, 2 \epsilon).
\end{align*}
Furthermore, by \eqref{fa4}, for any $0 \leq s < l$,
\begin{align*}
&\sum_{t=0}^{m\left(c_{s}(y)\right)-1}\varphi\left(\pi \circ f_{w}^{t_s(y)+1,t}\left(x_{s}\right)\right)\\
&=\sum_{t=0}^{m(c_{s}(y))-1}\varphi(g_{w}^{1,t_s(y)+t}(y))+\sum_{t=0}^{m(c_{s}(y))-1}\left(\varphi(\pi\circ f_{w}^{t_s(y)+1,t}(x_s))-\varphi(\pi\circ f_{w}^{1,t_s(y)+t}(x))\right)+ \\
&\sum_{t=0}^{m(c_{s}(y))-1}\left(\varphi(\pi\circ f_{w}^{1,t_s(y)+t}(x))-\varphi(g_{w}^{1,t_s(y)+t}(y))\right) \\
&=\sum_{t=0}^{m(c_{s}(y))-1}\varphi(g_{w}^{1,t_s(y)+t}(y))+\sum_{t=0}^{m(c_{s}(y))-1}\left(\varphi(\pi\circ f_{w}^{t_s(y)+1,t}(x_s))-\varphi(\pi\circ f_{w}^{1,t_s(y)+t}(x))\right)+ \\
&\sum_{t=0}^{m(c_{s}(y))-1}\left(\varphi(g_{w}^{1,t_s(y)+t}\circ\pi(x))-\varphi(g_{w}^{1,t_s(y)+t}(y))\right)\\
& \leq \sum_{t=0}^{m(c_{s}(y))-1} \varphi\left(g_{w}^{1,t_s(y)+t}(y)\right) + m(c_{s}(y)) \operatorname{var}(\varphi \circ \pi, 2 \epsilon) + m(c_{s}(y)) \epsilon.
\end{align*}
Similarly, when $ s=l $, we have
\[
\sum_{t=0}^{n-t_l(y)} \varphi\left(\pi \circ f_{w}^{t_l(y)+1,t}(x_l)\right)\leq \sum_{t=0}^{n-t_l(y)} \varphi\left(g_{w}^{1,t_l(y)+t}(y)\right) + (n+1-t_l(y))\left(\operatorname{var}(\varphi \circ \pi, 2 \epsilon) + \epsilon\right).
\]
Hence, 
\begin{align*}
& S_{w, n} \varphi \circ \pi(x) \\
\leq &\sum_{s = 0}^{l-1}\left(\sum_{t=0}^{m(c_{s}(y))-1} \varphi\left(g_{w}^{1,t_s(y)+t}(y)\right) + m(c_{s}(y)) \operatorname{var}(\varphi \circ \pi, 2 \epsilon) + m(c_{s}(y)) \epsilon\right)+\\
& \sum_{t=0}^{n-t_l(y)} \varphi\left(g_{w}^{1,t_l(y)+t}(y)\right) + (n+1-t_l(y))\left(\operatorname{var}(\varphi \circ \pi, 2 \epsilon) + \epsilon\right)+ (n+1) \operatorname{var}(\varphi \circ \pi, 2 \epsilon).
\end{align*}
That is
\begin{equation}\label{fa7}
S_{w, n} \varphi \circ \pi(x) \leq S_{w, n} \varphi(y) + 2(n+1) \operatorname{var}(\varphi \circ \pi, 2 \epsilon) + (n+1) \epsilon.
\end{equation}
Hence, using \ref{b} and \eqref{fa7},
\[
\sum_{x \in H_{w, n}} e^{S_{w, n} \varphi \circ \pi(x)} \leq \sum_{y \in E_{w, n}} \#\left(\mathcal{V}_{w,y}\right) \exp \left(S_{w, n} \varphi(y) + 2(n+1) \operatorname{var}(\varphi \circ \pi, 2 \epsilon) + (n+1) \epsilon\right). 
\]
Then,
\begin{align*}
& \quad P_n(\Phi;w, \varphi \circ \pi,4\epsilon)=\sup_E\left\{\sum_{x\in E}e^{S_{w,n}\varphi(x)}:E\text{ is an }(n,w,4\epsilon;\Phi)\text{-separated set for }X\right\} \\
& \leq \sum_{y \in E_{w, n}} \#\left(\mathcal{V}_{w,y}\right) \exp \left(S_{w, n} \varphi(y) + 2(n+1) \operatorname{var}(\varphi \circ \pi, 2 \epsilon) + (n+1) \epsilon\right).
\end{align*}
Thus,
\begin{align*}
& \quad \frac{1}{\#(I^{1,n})} \sum_{w \in I^{1,n}} P_n(\Phi;w, \varphi \circ \pi,4\epsilon)\\
& \leq \frac{1}{\#(I^{1,n})} \sum_{w \in I^{1,n}} \sum_{y \in E_{w, n}} \#\left(\mathcal{V}_{w,y}\right) \exp \left(S_{w, n} \varphi(y) + 2(n+1) \operatorname{var}(\varphi \circ \pi, 2 \epsilon) + (n+1) \epsilon\right).
\end{align*}
From \eqref{fa6}, 
\begin{align*}
&\quad P_{n}(\Phi;\varphi \circ \pi,  4 \epsilon) \\
\leq & \frac{1}{\#(I^{1,n})} \sum_{w \in I^{1,n}} \sum_{y \in E_{w, n}} e^{(n+M)(a + \alpha)} \exp \left(S_{w, n} \varphi(y) + 2(n+1) \operatorname{var}(\varphi \circ \pi, 2 \epsilon) + (n+1) \epsilon\right)\\
\leq & \exp\left((n+M)(a + \alpha) + 2(n+1) \operatorname{var}(\varphi \circ \pi, 2 \epsilon) + (n+1) \epsilon \right) \cdot \frac{1}{\#(I^{1,n})} \sum_{w \in I^{1,n}} \sum_{y \in E_{w, n}} e^{S_{w, n} \varphi(y)}.
\end{align*}
Therefore, 
\begin{align*}
& \quad \log P_{n}(\Phi;\varphi \circ \pi,  4 \epsilon) \\
& \leq (n+M)(a + \alpha) + 2(n+1) \operatorname{var}(\varphi \circ \pi, 2 \epsilon) + (n+1) \epsilon + \log \frac{1}{\#(I^{1,n})} \sum_{w \in I^{1,n}} \sum_{y \in E_{w, n}} e^{S_{w, n} \varphi(y)}.
\end{align*}
And then connecting with \eqref{fa3}, we have 
\begin{align*}
\frac{1}{n} \log P_{n}(\Phi;\varphi \circ \pi,  4 \epsilon)
 \leq & \frac{n+M}{n}(a + \alpha)+ 2\frac{n+1}{n} \operatorname{var}(\varphi \circ \pi, 2 \epsilon) + \frac{n+1}{n}\epsilon+ \\
& \frac{1}{n} \log \frac{1}{\#(I^{1,n})} \sum_{w \in I^{1,n}} \sum_{y \in E_{w, n}} e^{S_{w, n} \varphi(y)} \\
< & \frac{n+M}{n}(a + \alpha) + 2\frac{n+1}{n} \operatorname{var}(\varphi \circ \pi, 2 \epsilon) + \frac{n+1}{n}\epsilon + P(\Psi, \varphi)+\beta.
\end{align*}
Let $n \to \infty$ and $\epsilon \to 0$, 
\[
P(\Phi, \varphi \circ \pi) \leq P(\Psi, \varphi) + a + \alpha+\beta.
\]
Since $\alpha$ and $\beta$ are arbitrary, this implies $ P(\Phi, \varphi \circ \pi) \leq P(\Psi, \varphi) + a $ and the theorem is proved.
\end{proof}
\end{theorem}
\begin{remark}
If $ \Phi^{(i)}=\mathcal{F}:=\{f_{0},f_{1},\ldots,f_{m-1}\} $ and $ \Psi^{(i)}=\mathcal{G}:=\{g_{0},g_{1},\ldots,g_{m-1}\} $ for all $ i \in \mathbb{N}$, then Theorem \ref{con} yields the following inequality for topological pressures with a factor map of free semigroup actions:
\[
P(\mathcal{F},\varphi\circ\pi)\leq P(\mathcal{G},\varphi)+\sup_{y\in Y}H(\mathcal{F};\pi^{-1}(y)), 
\]
In comparison, the main result in \cite{Zhao} gives
\[
P(\mathcal{F},\varphi\circ\pi)\leq P(\mathcal{G},\varphi)+\sup_{y\in Y}h_{\text{Bi\'s}}(\mathcal{F};\pi^{-1}(y)), 
\]
where $h_{\text{Bi\'s}}(\mathcal{F};K)$ denotes the Bi\'s entropy of a compact subset $K$ of $X$. Since it is straightforward to verify that
$H(\mathcal{F};\pi^{-1}(y)) \leq h_{\text{Bi\'s}}(\mathcal{F};\pi^{-1}(y))$ for every $y \in Y$. Therefore, our result constitutes a refinement and generalization of the inequality in \cite{Zhao}. 
\end{remark}

%
%

%
%
%

\end{document}